\documentclass{article}
\usepackage[T1]{fontenc}
\usepackage{stackengine}
  \usepackage{setspace}
\usepackage{IEEEtrantools}
\usepackage{cite, authblk}
\usepackage{amssymb,amsfonts, anyfontsize, soul}
\usepackage{enumitem}

\usepackage{amsmath, amsthm}               

  {

  }
\IEEEeqnarraydefcolsep{0}{\leftmargini} 

\usepackage{algorithmic}
\usepackage{graphicx}
\usepackage{textcomp}
\usepackage{xcolor}

\usepackage{algorithm,algorithmic, longtable, xtab}
\usepackage{mathtools, stackengine, mathrsfs, setspace, bm, multirow, tikz, pgfplots, adjustbox, calc, url}
\usepackage[colorlinks=true, citecolor=blue, linkcolor=black!20!red]{hyperref}

\usepackage{xparse}

\makeatletter

\let\oldlt\longtable
\let\endoldlt\endlongtable
\def\longtable{\@ifnextchar[\longtable@i \longtable@ii}
\def\longtable@i[#1]{\begin{figure}[t]
\onecolumn
\begin{minipage}{0.5\textwidth}
\oldlt[#1]
}
\def\longtable@ii{\begin{figure}[t]
\onecolumn
\begin{minipage}{0.5\textwidth}
\oldlt
}
\def\endlongtable{\endoldlt
\end{minipage}
\twocolumn
\end{figure}}
\makeatother
\allowdisplaybreaks

\begin{document}
\title{Flexible Ramping Product Procurement \\in Day-Ahead Markets}
%
%
%

\author{Ogun~Yurdakul$^{1}$, {Erik Ela}$^{2}$,~Sahin~Albayrak$^{1}$\\
{\fontsize{9}{5}\selectfont
$^{1}$Department of Electrical Engineering and Computer Science, Technical University of Berlin, Berlin, Germany (e-mail: yurdakul, sahin.albayrak@tu-berlin.de). \newline
$^{2}$Power Delivery and Utilization Department, Electric Power Research Institute, Palo Alto, CA 94304-1395 USA (e-mail: eela@epri.com).}
}
\providecommand{\keywords}[1]{\textbf{\textit{Index terms---}} #1}

\date{}
\maketitle

\begin{abstract}
This article puts forward a methodology for procuring flexible ramping products ({FRP}s) in the day-ahead market ({DAM}). The proposed methodology comprises two market passes, the first of which employs a stochastic unit commitment ({SUC}) model that explicitly evaluates the uncertainty and the intra-hourly and inter-hourly variability of net load so as to minimize the expected total operating cost. The second pass clears the {DAM} while imposing {FRP} requirements. The cornerstone of our work is to set the {FRP} requirements at levels that drive the {DAM} decisions toward the optimal {SUC} decisions. Our methodology provides an economic underpinning for the stipulated {FRP} requirements, and it brings forth {DAM} awards that reduce the costs toward the expected total cost under {SUC}, while conforming to the chief {DAM} design principles. By preemptively considering the dispatch costs before awarding {FRP}s, it can further avert unexpectedly high costs that could result with the deployment of procured {FRP}s. We conduct numerical studies and lay out the relative merits of the proposed methodology vis-à-vis selected benchmarks based on various evaluation metrics. 
\end{abstract}
\newpage
\section*{Nomenclature}\label{0}
\subsection*{Sets,Indices}\label{1}
\begin{tabular}[h]{ll}
$\mathscr{H}, h$  & \hspace{0.45cm}set, index of hourly periods\\
$\mathscr{K}, k$  & \hspace{0.45cm}set, index of intra-hourly subperiods\\
$\mathscr{N}, n$  & \hspace{0.45cm}set, index of nodes\\
$\mathscr{L}, \ell$  & \hspace{0.45cm}set, index of lines\\
$\mathscr{G}, g$  & \hspace{0.45cm}set, index of dispatchable generators ({DG}s)\\
$\mathscr{S}_g, s$ & \hspace{0.45cm}set, index of piecewise cost intervals for {DG} $g$
\end{tabular} 
\subsection*{Variables}\label{2}
\begin{tabular}[h]{ll}
$u_{g}$ & \hspace{0.09cm} $\in \{0,1\}$, commitment status of {DG} $g$ \\
$v_{g}$ & \hspace{0.09cm} $\in \{0,1\}$, startup status of {DG} $g$ \\
$w_{g}$ & \hspace{0.09cm} $\in \{0,1\}$, shutdown status of {DG} $g$ \\
$p_{g}$ & \hspace{0.09cm} power generated above minimum by {DG} $g$ ({MW}) \\
$p_g^{s}$ & \hspace{0.09cm} power from linear segment $s$ for {DG} $g$ ({MW})\\
$p_{c}^{n}$& \hspace{0.09cm} curtailed load at node $n$ ({MW})\\
${d}_{n}$ & \hspace{0.09cm} net demand at node $n$ ({MW})\\
$r_{g}^{\uparrow}$& \hspace{0.09cm} upward flexible ramping product (up-{FRP}) provided  \\
& \hspace{0.09cm} by {DG} $g$ ({MW})\\
$r_{g}^{\downarrow}$& \hspace{0.09cm} downward {FRP} (down-{FRP}) provided by {DG} $g$ ({MW})\\
$r^{\uparrow}_{\wr} (r^{\downarrow}_{\wr})$& \hspace{0.09cm} up-{FRP} (down-{FRP}) shortfall ({MW})\\
$\lambda^{n}$ & \hspace{0.2cm}locational marginal price at node $n$ (dual) ({\$/MWh}) \\
$\varphi^{\uparrow}(\varphi^{\downarrow})$ & \hspace{0.2cm}up-{FRP} (down-{FRP}) price (dual) ({\$/MWh})\\
\end{tabular} 
\subsection*{Parameters}\label{3}
\begin{tabular}[h]{ll}
$p_{g}^{\circ}$ & \hspace{-0.25cm}initial power generated above minimum by {DG} $g$ ({MW})\\
$u_{g}^{\circ}$ & \hspace{-0.25cm}$\in \{0,1\}$, initial commitment status of {DG} $g$\\
$T^{\uparrow}_{g}(T^{\downarrow}_{g})$ & \hspace{-0.25cm}minimum uptime (downtime) of {DG} $g$ (h)\\
$T^{\uparrow, \circ}_{g}(T^{\downarrow, \circ}_{g})$ & \hspace{-0.25cm}number of hours {DG} $g$ has been online (offline) \\
& \hspace{-0.25cm}before the scheduling horizon (h)\\
$\overline{P}_{g}^{s}$ & \hspace{-0.25cm}maximum power available from piecewise segment\\
& \hspace{-0.25cm}$s$ for {DG} $g$ ({MW})\\
$\overline{P}_{g}(\underline{P}_{g})$ & \hspace{-0.25cm}maximum (minimum) power output of {DG} $g$ ({MW})\\
$\alpha_{g}^{s}$ & \hspace{-0.25cm}cost coefficient for piecewise segment $s$ for \\
& \hspace{-0.25cm}{DG} $g$ ({\$/MWh})\\
\end{tabular} 
\newpage
\begin{tabular}[h]{ll}
\hspace{-0.81cm} $\alpha^{c}$ & \hspace{-0.21cm}penalty cost for load curtailment ({\$/MWh})\\
\hspace{-0.81cm} $\alpha^{r}$ & \hspace{-0.21cm}penalty cost for {FRP} shortfall ({\$/MWh})\\
\hspace{-0.81cm} $\alpha_{g}^{v}$ & \hspace{-0.21cm}startup cost of {DG} $g$ (\textit{\$})\\
\hspace{-0.81cm} $\alpha_{g}^{u}$ & \hspace{-0.21cm}cost of {DG} $g$ running and operating at $\underline{P}_{g}$ ({\$/h})\\
\hspace{-0.81cm} $\Delta_{g}^{\uparrow} (\Delta_{g}^{\downarrow})$ & \hspace{-0.21cm}ramp-up (ramp-down) rate limit of {DG} $g$ ({MW/h})\\
\hspace{-0.81cm} $\Delta_{g}^{\uparrow, \circ} (\Delta_{g}^{\downarrow, \circ})$ & \hspace{-0.21cm}startup (shutdown) rate limit of {DG} $g$ ({MW})\\
\hspace{-0.81cm} $\overline{f}_{\ell}(\underline{f}_{\ell})$ & \hspace{-0.21cm}maximum (minimum) real power flow allowed\\
\hspace{-0.81cm} & \hspace{-0.21cm}on line $\ell$ ({MW})\\
\hspace{-0.81cm} $\Psi_{\ell}^{n}$ & \hspace{-0.21cm}injection shift factor of line $\ell$ with respect to node $n$ \\
\hspace{-0.81cm} & \hspace{-0.21cm}by {DG} $g$ ({MW})\\
\hspace{-0.81cm} $\hat{d}_{n}$&\hspace{-0.21cm}net demand at node $n$ ({MW})\\ 
\hspace{-0.81cm} $\xi^{n}_{i}$ &\hspace{-0.21cm}net load for realization $i$ at node $n$ ({MW})\\
\hspace{-0.81cm} $\underline{\rho}^{\uparrow}(\underline{\rho}^{\downarrow})$ & \hspace{-0.21cm}system-wide up-{FRP} (down-{FRP}) requirement ({MW})
\end{tabular} 
\section{Introduction}\label{1}
\IEEEPARstart{W}{ith} the deepening penetration of variable energy resources ({VER}s) into electric power systems, the uncertainty in net load, i.e., system load less {VER} generation, is expected to further increase. Since neither system load nor {VER} generation levels can be known with perfect accuracy ahead of time, the net load levels considered by day-ahead planning tools are inherently uncertain, often materializing differently than predicted. \par
{VER}s further exacerbate the variability in net load, which refers to changing net load levels that are anticipated. We bifurcate net load variability based on its temporal granularity to intra-hourly and inter-hourly net load variability. Since day-ahead markets ({DAM}s) and residual (or reliability) unit commitment ({RUC}) processes are executed based on hourly snapshots, whereas real-time markets ({RTM}s) are cleared with a finer temporal granularity, the intra-hourly variation of net load is not recognized by most day-ahead planning tools. On the flip side, {DAM}s inherently take into account the inter-hourly variation of net load by virtue of solving a time-coupled multi-interval unit commitment ({UC}) problem. Nevertheless, since system operators\footnote{We use the generic term \textit{system operator} to represent the independent entity responsible for operating and controlling the transmission network---be it an independent system operator ({ISO}), regional transmission organization ({RTO}), or a transmission system operator.} ({SO}s) use the bid-in, not the forecast, {net} demand to clear the {DAM}, the {DAM} energy schedules may not be able to meet steep inter-hourly net load variations.\par
Repercussions of the uncertainty and variability of net load include intervals with a ramp shortage, during which system assets lack the ramping capability for meeting the load in the subsequent interval. Such ramp shortages may force {SO}s to use up spinning or regulation reserves, thus jeopardizing power system reliability, and may also yield price spikes in {RTM}s. Previous studies maintain that such transient extreme prices are not necessarily indicative of an underlying chronic scarcity of ramping capability, but rather point to the shortcomings of the utilized scheduling methods \cite{ramp:malley}. To avoid ramp shortages, {SO}s may make out-of-market ({OOM}) adjustments, such as withholding some generation capacity, 
{committing} 
{dispatchable generators (DGs)} outside of market mechanisms, and the manual modification of load forecast. Resorting to {OOM} actions to engage with ramp shortages not only cloaks the severity of the problem, but also undermines the development of an efficient market-based solution that should work hand and glove with existing market mechanisms, {including providing the necessary incentives that may attract existing and emerging resources to offer flexibly}. \par
In recent years, California ISO ({CAISO}) and Midcontinent ISO ({MISO}) augmented their markets with ``flexiramp'' \cite{ramp:caisofrpref} and ``ramp capability'' \cite{ramp:miso} products, respectively, which is a {short-term rampable} capacity designated to meet a level of net load higher or lower than expected in subsequent time periods. We henceforth refer to such products for upward and downward ramping capability by upward flexible ramping product (up-{FRP}) and downward {FRP} (down-{FRP}), respectively. Separate 
{non-zero availability} bids are not submitted for {FRP}s, and it assumed that a resource submitting an energy offer is willing to provide both energy and {FRP}, unless it chooses to opt out of clearing for {FRP}s. The market clearing model co-optimizes energy and {FRP} quantity with other ancillary services and sets the price for {FRP}s based on the lost opportunity cost of the marginal unit, that is, the foregone net profit from the energy market that the resource could have earned, had it not received an {FRP} award.\par
The effective acquisition of {FRP}s calls for approaches that set {FRP} requirements so as to safeguard against the uncertainty and the inter-hourly and intra-hourly variability of net load. As shall be described in Section \ref{2}, many methods procure {FRP}s based on the 95\% confidence interval of net load. These methods set rules that influence the decisions of an optimization problem that is inherently stochastic, without assessing how such rules will impact the expected cost function. Absent economic underpinnings, imposing rules in markets that lay down specific numerical thresholds may lack a warranted basis and so seem arbitrary. Regulators and {SO}s need to have a clear picture of the costs and benefits of FRPs to ratepayers so that they can render proper value judgments on FRP requirements. 
 \par
Another major drawback of most {FRP} procurement methodologies is the lack of consideration for the production costs that would be incurred, were the procured {FRP}s to be deployed. Since {FRP}s are expected to be called on frequently, methodologies need to be devised that can take account of the production costs while awarding {FRP}s. Further, these methodologies must be able to be brought into force with few modifications to existing market mechanisms and send transparent price signals that provide incentives for the provision of {FRP}s when required and ultimately give impetus to investments in {DG}s with a fast-ramping capability. For instance, one of the approaches considered as part of the Day-Ahead Market Enhancements ({DAME}) program of {CAISO} was the joint utilization of bid-in demand and forecast load in a single market pass to clear energy, {FRP}, and other ancillary services, which faced criticism from stakeholders for being overly reliant on the accuracy of load forecasts in setting the prices \cite{ramp:caisonew}.\par
In this paper, we map out a methodology for {FRP} acquisition in the {DAM}, which comprises two market passes that are executed sequentially. The first pass is a stochastic unit commitment ($\mathsf{SUC}$) problem that explicitly models the uncertainty in net load and is formulated at an intra-hourly granularity. While the optimal solution to the $\mathsf{SUC}$ problem would theoretically yield the minimum expected total cost, both its stochastic nature and its temporal granularity do not conform to existing {DAM} designs. \par
To surmount these challenges in the design of our methodology, we marry the quantitative insights obtained from the optimal $\mathsf{SUC}$ solution with today’s {DAM} models. We evaluate the {FRP} requirements based on the optimal $\mathsf{SUC}$ solution and impose them in the second pass, which is a {DAM} clearing ($\mathsf{DAMC}$) problem with little modification to current {DAM} models. Our design choices impart the capability to procure {FRP}s in a way that brings {DAM} schedules toward the optimal $\mathsf{SUC}$ solution, thereby reducing the expected total cost. \par
We draw upon the optimal $\mathsf{DAMC}$ decisions to clear the {DAM} and set the prices: we harness the shadow price of the nodal power balance constraint to set the locational marginal price ({LMP}) of energy and the dual variable of the up and down-{FRP} procurement constraints for pricing up and down-{FRP}, respectively. By utilizing the bid-in demand in the second pass, we preclude the energy and {FRP} prices from being directly exposed to the {SO} forecasts. Owing to the explicit incorporation of transmission constraints in both passes, we further provide for the nodal deliverability of procured {FRP}s for all considered scenarios.
\subsection{Contributions and Structure of the Paper}
Broadly, the contributions of this paper are as follows:
\begin{enumerate}
 \item We work out a methodology for {FRP} procurement that sets the {FRP} requirements and clears the {DAM} so as to drive the {DAM} decisions toward the optimal decisions of the $\mathsf{SUC}$ problem, which assesses the uncertainty and the inter-hourly and intra-hourly variability of net load. As such, the stipulated {FRP} requirements are explicitly buttressed by an economic pillar. 
\item In awarding {FRP}s, the proposed approach is equipped to take explicit consideration of the dispatch costs that will arise due to the deployment of the procured {FRP}s. 
\item The proposed methodology adheres to the key principles of the {DAM}: it clears {FRP}s and energy schedules at an hourly granularity and leverages the optimal dual variables to set the {FRP} and energy prices. 
\end{enumerate}
The remainder of this paper contains four sections. In Section \ref{2}, we provide an overview of {FRP} practices and designs from the industry and academic literature. In Section \ref{3}, we lay out the analytical underpinnings of our methodology and present our problem formulations. We describe the numerical studies conducted to demonstrate the application and effectiveness of the proposed methodology in Section \ref{4} and discuss the results. We summarize the paper and provide directions for future studies in Section \ref{5}. 
\section{{FRP}s: Background and Related Work}\label{2}
In this section, we 
provide an overview of the approaches to {FRP} design and procurement that are proposed in the academic literature and implemented by {SO}s.
\subsection{Literature Review}\label{2b}
In recent years, {FRP}s have received considerable attention in the literature (see \cite{ramp:survey} and the references therein for a comprehensive review). The studies carried out in \cite{ramp:hobbs,ramp:malley,ramp:rosenwald,ramp:spin} focus on {FRP} procurement in {RTM}s, which is pivotal especially because the {RTM} may clear based on a single-interval dispatch model, or even when it clears based on a multi-interval dispatch model, the dispatch and prices for only one interval may be binding, with those for the remaining intervals acting as advisory. \par
In \cite{ramp:hobbs}, the authors set up a real-time {UC} formulation that procures up and down-{FRP}s based on the forecast change in net load and an assumed error range, and they set the {FRP} prices by using the dual variable of the {FRP} procurement constraints. They compare their approach with a multi-stage stochastic {UC} problem and show that their approach may underestimate the cost of relying on an expensive generator for {FRP} provision, when it disregards the dispatch costs that would arise should an expensive generator receiving an {FRP} award be called on to generate power. The study in \cite{ramp:malley} aims to develop an incentive-compatible pricing method in a multi-interval dispatch model for {RTM}s to meet the variability of net load but does not address net load uncertainty. \par
In contrast to \cite{ramp:hobbs,ramp:malley}, \cite{ramp:rosenwald} incorporates a demand curve for {FRP} into its problem formulation, considers {FRP} procurement on a both zonal and system-wide basis, and studies the influence of having {FRP} availability offers on cleared {FRP} quantities and prices. Reference \cite{ramp:spin} lays out a methodology to use a portion of the spinning reserve committed in the {DAM} for {FRP} procurement in the {RTM}, which assesses the net load variation between successive intervals jointly with the 95\% confidence interval for net load forecast error.	\par
References \cite{ramp:mojdeh, ramp:env} address {FRP} procurement in the {DAM}, for which the studies need to take into account the granularity difference with the {RTM} and recognize that the cleared demand may deviate from the load forecast of the {SO}. To this end, the approach proposed in \cite{ramp:mojdeh} augments the hourly {FRP} requirements to ensure sufficient ramp capability for each of the four 15-minute intervals within the hour. Nevertheless, \cite{ramp:mojdeh} does not distinguish between bid-in and forecast load, and it assumes that the net load forecast errors follow a pre-specified distribution and makes use of the 95\% confidence interval in acquiring {FRP}s. Reference \cite{ramp:env} takes as its primary focus the environmental ramifications of {FRP} procurement and leverages the dual variables of the {FRP} procurement constraints to set the {FRP} prices, as adopted in \cite{ramp:hobbs,ramp:malley,ramp:rosenwald,ramp:spin,ramp:mojdeh}. However, \cite{ramp:env} fails to distinguish between bid-in and forecast load and, similar to\cite{ramp:spin,ramp:mojdeh}, utilizes the 95\% confidence interval to acquire {FRP}s for meeting net load uncertainty. 
\subsection{Industry Practices}\label{2c}
The incorporation of {FRP}s into wholesale electricity markets was spearheaded by {CAISO} and {MISO} in 2010s. Currently, {CAISO} acquires {FRP}s only in its {RTM}, for which it does not accept separate financial bids and relies on a demand curve that considers both the variability and uncertainty of net load \cite{ramp:caisofrp}. For incorporating {FRP}s into the {DAM}, discussions among {CAISO} and other stakeholders have been underway as part of the {DAME} program. In the published {DAME} straw proposals\cite{ramp:caisonew}, {CAISO} put forth the integration of reliability capacity up/down products to the {DAM} for meeting the ramping requirement difference between the cleared and forecast hourly net load. {CAISO} further proposes the integration of imbalance reserves up/down products to the {DAM} so as to meet the imbalances due to net load uncertainty and the granularity difference between hourly and 15-minute forecasts.\par
In contrast, {MISO} procures {FRP}s in both the {DAM} and {RTM}. In its market clearing model, {MISO} employs a separate single-segment demand curve for up and down-{FRP} and co-optimizes {FRP}s with energy and other ancillary services \cite{ramp:miso}. By assessing both the variability and uncertainty of net load, {MISO} determines {FRP} requirements and leverages the dual variable of the {FRP} procurement constraints to set the price for {FRP}s. {In addition, Southwest Power Pool recently introduced FRPs to its RTM, which are designed to manage the short-term variability and uncertainty in net load} \cite{ramp:spp}. \par
{Note that while the FRPs currently in force are predominantly geared to meet the short-term variability and uncertainty in net load, the discussion in Section {\ref{1}} made clear that there is a substantial need to explore how similar products can be designed and leveraged so as to engage with the challenges due to said variability and uncertainty in the day-ahead time frame. In the next section, we map out our methodology for FRP procurement in the day-ahead time frame, which aims to address precisely these challenges.}
\section{Analytical Underpinnings of the\\ Proposed Methodology}\label{3}
In this section, we set forth the analytical underpinnings of the proposed {FRP} procurement methodology. 
We discretize the time axis and---commensurate with the typical planning horizon of {DAM}s---adopt 24 hours as the scheduling horizon. We denote by $h$ the index for the hourly periods of a day and construct the set $\mathscr{H} \coloneqq \{ h \colon h=1,\ldots,24\}$. We consider $K$ intra-hourly subperiods of equal duration for each hour $h$ and introduce the index $k$ for intra-hourly subperiods, such that $k=1,\ldots,24\cdot K$. We denote by $\zeta$ the duration of each intra-hourly subperiod in minutes. For each hour $h\in \mathscr{H}$, we construct the set $\mathscr{K}_h \coloneqq \{ k \colon k=(h-1)\cdot K + 1,\ldots,h\cdot K\}$. We define the set $\mathscr{K}\coloneqq\cup_{h\in\mathscr{H}} \mathscr{K}_h$. Note that both the set $\mathscr{H}$ and the set $\mathscr{K}$ represent time periods that span over 24 hours, yet their cardinalities are different, as their representations adopt different temporal granularities. While the set $\mathscr{H}$ contains 24 elements representing 24 hourly periods, the set $\mathscr{K}$ contains $24\cdot K$ elements for $24\cdot K$ intra-hourly time periods, each with a $\zeta$-minute duration. \par
Since the problem formulations that will be laid out in the next subsections have different time characteristics, we distinguish the temporal granularity of variables in our notation. We denote by $x[k]$ the variable $x$ in the intra-hourly subperiod $k$, which adopts $\zeta$ minutes as the smallest indecomposable unit of time and assumes that system conditions remain constant over each subperiod $k$. On the flip side, the term $x(h)$ denotes the variable $x$ in hour $h$, which has an hourly granularity, thus assuming that system conditions hold constant over each hour $h \in \mathscr{H}$. 
The primary objective in the design of the proposed methodology is to explicitly assess the expected costs under the uncertainty and variability of net load---all the while observing the chief design principles of {DAM}s. To achieve this objective, we design our methodology such that it consists of two market passes that are executed back to back. \par
The first pass is the ($\mathsf{SUC}$) problem, which expressly recognizes the uncertainty in net load and is solved at an intra-hourly temporal granularity. Optimal solutions to the $\mathsf{SUC}$ problem are leveraged in setting the {FRP} requirements in the second pass $\mathsf{DAMC}$ problem, which is a deterministic {UC} problem solved to clear the {DAM}, evaluate energy and {FRP} awards, and set the up-{FRP} and down-{FRP} prices. We next describe each market pass in turn.
\subsection{Stochastic Unit Commitment ($\mathsf{SUC}$) Problem Formulation}\label{3a}
In this subsection, we present the mathematical formulation of the $\mathsf{SUC}$ problem, which is modeled at a $\zeta$-minute granularity so as to assess both the inter-hourly and intra-hourly variability of net load. \par
The $\mathsf{SUC}$ problem explicitly models the uncertainty in net load. We represent the uncertain net load in the transmission network over the study period by ${\xi}$ defined on the probability space $(\Omega, \mathcal{F}, \mathcal{P})$, where $\Omega$ is a sample space, $\mathcal{F}$ is a set of subsets of $\Omega$ that is a $\sigma-$algebra, and $\mathcal{P}$ is a probability distribution on $\mathcal{F}$. Suppose we draw an independent and identically distributed (\textit{i.i.d.}) random sample $\{{\xi}_{1},\ldots,{\xi}_{i},\ldots,{\xi}_{I}\}$ of size $I$ of the random net load ${\xi}$. For ease of notation, we define the set $\mathscr{I} \coloneqq \{i \colon i = 1,\ldots,I\}$ and express by ${\xi}^{n}_{i}[k]$ the net load at node $n \in \mathscr{N}$ in the intra-hourly subperiod $k \in \mathscr{K}$ for the net load realization with index $i \in \mathscr{I}$. \par
The $\mathsf{SUC}$ problem is a two-stage mixed-integer program that reflects the sequence in which the {UC} and economic dispatch decisions are rolled out. We make use of the sample average approximation to express the first-stage problem as
\begin{IEEEeqnarray}{ll}
\hspace{-0.45cm}\underset{u_g[k], v_g[k]}{\text{minimize}} &\hspace{0.5cm}  \sum_{k \in \mathscr{K}} \sum_{g \in \mathscr{G}} \Big[ \alpha^{u}_{g} u_{g}[k] + \alpha^{v}_{g} v_{g}[k] \Big]  +  \frac{1}{I}\sum_{i \in \mathscr{I}}\mathcal{Q}({x},{\xi_{i}}), \label{obj}\\
\hspace{-0.45cm}\text{ subject to} & \nonumber
\end{IEEEeqnarray}

\begin{IEEEeqnarray}{lcl}
\hspace{-.49cm}u_{g}[k]-u_{g}[k-1] &\hspace{0.2cm}&\nonumber\\
\hspace{-.49cm}\hspace{1.5cm}= v_{g}[k] - w_{g}[k]&&\forall k \in \mathscr{K}\setminus\{1\}, \label{st}\\
\hspace{-.49cm}u_{g}[k]-u_{g}^{\circ} = v_{g}[k] - w_{g}[k] &&\forall k \in \{1\}, \label{st1}\\
\hspace{-.49cm}\sum_{k'=k-(T^{\uparrow}_{g}\cdot K)+1}^{k} \hspace{-0.5cm}v_{g}[k'] \leq u_{g}[k] &&\forall k \in \{T^{\uparrow}_{g} \cdot K,\ldots,24\cdot K\},\label{upt}  \\
\hspace{-.49cm}\sum_{k'=k-(T^{\downarrow}_{g}\cdot K)+1}^{k} \hspace{-0.5cm}w_{g}[k'] \leq1- u_{g}[k]&&\forall k \in \{T^{\downarrow}_{g}\cdot K,\ldots,24\cdot K\}, \label{dwt}  \\
\hspace{-.49cm}\sum_{k'=1}^{\min\{u_g^{\circ}\cdot(T^{\uparrow}_{g} - T^{\uparrow, \circ}_{g})\cdot K, \,24\cdot K\}} \hspace{-0.5cm}w_{g}[k'] = 0,&& \label{iupt}  \\
\hspace{-.49cm}\sum_{k'=1}^{\min\{(1-u_g^{\circ})\cdot(T^{\downarrow}_{g} - T^{\downarrow, \circ}_{g})\cdot K, \,24\cdot K\}} \hspace{-1.0cm}v_{g}[k'] = 0,&& \label{idwt}  \\
\hspace{-.49cm}u_{g}[k], v_{g}[k], w_{g}[k]  \in \{0,1\} && \forall k \in \mathscr{K} ,\label{bin}\\
\hspace{-.49cm}u_{g}[k] = u_{g}[h\cdot K]&&\forall k \in \mathscr{K}_h\setminus \{K\cdot h\}, \forall h \in \mathscr{H}, \label{int}
\end{IEEEeqnarray}
where \eqref{st}--\eqref{int} hold {for all {DG}s} $g \in \mathscr{G}$. The objective  \eqref{obj} of the first stage is to minimize the commitment and startup costs plus the expected dispatch and load curtailment costs. 
\par
We adopt a three-binary formulation for {DG}s 
and enforce by (\ref{st},~\ref{st1}) logical constraints that relate the variables $u[k]$, $v[k]$, and $w[k]$. The minimum uptime requirement of the the {DG}s is imposed by (\ref{upt},~\ref{iupt}), and the minimum downtime requirements of the the {DG}s is stated in (\ref{dwt},~\ref{idwt}). We enforce by \eqref{int} the constraint that the commitment variables attain equal values for all intra-hourly subperiods $k \in \mathscr{K}_h$ within each hour $h \in \mathscr{H}$. As such, \eqref{int} ensures that optimal $\mathsf{SUC}$ decisions can be implemented in the second pass by the $\mathsf{DAMC}$ problem, which, as will be spelled out in Section \ref{3b}, has an hourly granularity. We represent all first-stage decision variables by the vector ${x}$, which comprises $u_g[k]$, $v_g[k]$, and $w_g[k]$. For a specific vector of first-stage decision variables $x$ and net load realization $\xi_{i}$, {the value function} $\mathcal{Q}({x},\xi_{i})$ is computed by solving the following second-stage problem: 
\begin{IEEEeqnarray}{lcl}
\hspace{-0.82cm}\underset{\begin{subarray}{c}p_g[k], p_g^{s}[k], \\p^{n}_c[k]\end{subarray}}{\text{minimize}} && \sum_{k \in \mathscr{K}}\bigg[ \sum_{g \in \mathscr{G}} \sum_{s \in \mathscr{S}_g} {\alpha}^{s}_{g}p^{s}_{g}[k] + \sum_{n \in \mathscr{N}} \alpha_{c}p^{n}_{c}[k]\bigg], \label{obj2}\\[5pt]
\hspace{-0.82cm}\text{subject to} & \nonumber
\end{IEEEeqnarray}
\begin{IEEEeqnarray}{{-0.01pt}l}
\,p_{g}[k] \leq (\overline{P}_{g}-\underline{P}_{g})u_{g}[k], \label{gl} \\
\,p_{g}[k] \leq p_{g}^{\circ} + \frac{\Delta^{\uparrow}_{g}}{K}u_{g}^{\circ} +(\Delta^{\uparrow,\circ}_{g}-\underline{P}_{g})v_{g}[k],\label{rul1} \\
\,p_{g}[k] \geq p_{g}^{\circ} - \frac{\Delta^{\downarrow}_{g}}{K}u_{g}^{\circ}-(\Delta^{\uparrow,\circ}_{g}-\frac{\Delta^{\uparrow}_{g}}{K}-\underline{P}_{g})w_{g}[k],\label{rll1} \\
\,p_{g}[k] \leq p_{g}[k-1] + \frac{\Delta^{\uparrow}_{g}}{K}u_{g}[k-1]+(\Delta^{\uparrow,\circ}_{g}-\underline{P}_{g})v_{g}[k],\label{rul} \\
\,p_{g}[k] \geq p_{g}[k-1] - \frac{\Delta^{\downarrow}_{g}}{K}u_{g}[k-1]-(\Delta^{\uparrow,\circ}_{g}-\frac{\Delta^{\uparrow}_{g}}{K}-\underline{P}_{g})w_{g}[k],\label{rll} \\
\,p_{g}[k] = \sum_{s \in \mathscr{S}_g} p_{g}^{s}[k] \label{ls} \\
\,p^{s}_{g}[k] \leq \overline{P}^{s}_{g} - \overline{P}^{\,s-1}_{g}\hspace{4.1cm} \forall {s} \in \mathscr{S}_g,\label{lsl} \\
\,p_g[k] \in \mathbb{R}_{+}, \label{nn1}\\
\,p^{s}_g[k] \in \mathbb{R}_{+}\hspace{5.25cm} \forall {s} \in \mathscr{S}_g , \label{nn2}\\
\, p^{n}[k] = \sum_{g \in \mathscr{G}^{n}} p_{g}[k] + p^{n}_{c}[k] - {\xi}_{i}^{n}[k]\hspace{2.15cm} \forall n \in \mathscr{N}, \label{netpower}\\
\, \sum_{n \in \mathscr{N}} p^{n}[k]  = 0 , \label{pb}\\
\, \underline{f}_{\ell} \leq \sum_{n \in \mathscr{N}}  \Psi_{\ell}^{n}p^{n}[k] \leq \overline{f}_{\ell}\hspace{3.5cm}\forall \ell \in \mathscr{L}, \label{pf}\\
\,p^{n}_c[k] \in \mathbb{R}_{+} \hspace{5.3cm}\forall n \in \mathscr{N},\label{nn3}
\end{IEEEeqnarray}
where \eqref{gl} and \eqref{ls}--\eqref{nn3} hold for all intra-hourly subperiods $k \in \mathscr{K}$, \eqref{rul1}--\eqref{rll1} for $k \in \{1\}$, \eqref{rul}--\eqref{rll} for $k \in \mathscr{K} \setminus\{1\}$, and \eqref{gl}--\eqref{nn2} for all {DG}s $g \in \mathscr{G}$. The second-stage problem \eqref{obj2}--\eqref{nn3} seeks to determine the power dispatch of the {DG}s ($p_g[k]$), power from each linear segment ($p_g^{s}[k]$) and the curtailed load ($p^{n}_c[k]$) for all $k \in \mathscr{K}$ with the objective \eqref{obj2} to minimize the dispatch costs of {DG}s and the penalty cost incurred due to load curtailment. For each net load realization $\xi_{i}$, we succinctly express by $y_i$ the second-stage decision variables, \textit{viz.}: $p_g[k]$, $p_g^{s}[k]$, and $p^{n}_c[k]$. We model the generator dispatch cost to be convex and piecewise linear in $p_g[k]$ in \eqref{obj2} based on the formulation presented in \cite{knuostwat}.\par
Equations \eqref{gl}--\eqref{rll} enforce {DG} generation and ramping limits. While \eqref{gl}--\eqref{rll1} are based on the formulation laid out in \cite{morlatram}, they additionally take into account the case that the startup (resp. shutdown) rate limit of a {DG} may be higher than its ramp-up (resp. ramp-down) rate limit and further explicitly model the initial subperiod of the simulation. We constrain the power output above minimum for each {DG} $g$, \textit{viz.}, $p_g[k]$, based on its maximum and minimum power output and commitment status in \eqref{gl}. We enforce that each {DG} respects its startup and ramp-up rate limits in the initial subperiod of the simulation and in the remaining subperiods of the simulation by \eqref{rul1} and \eqref{rul}, respectively. Similarly, \eqref{rll1} and \eqref{rll} ensure for each {DG} $g$ that the power output above minimum is within its shutdown and ramp-down rate limits in the initial subperiod of the simulation and in the remaining subperiods of the simulation, respectively. For each {DG} $g$, \eqref{ls} defines the relationship between power output above minimum and the power from each linear segment $s \in \mathscr{S}_g$, \textit{viz.}, $p_g^{s}[k]$, and \eqref{lsl} constrains $p_g^{s}[k]$ by the maximum power available from each linear segment, \textit{viz.}, $\overline{P}^{s}_{g}$, for $s \in \mathscr{S}_g$.\par
We state in \eqref{netpower} the net real power injection at each node $n \in \mathscr{N}$ with the convention that $p^{n}[k]>0$ if real power is injected into the system and lay down the system-wide power balance constraint in \eqref{pb}. We adopt the DC power flow model to represent the transmission constraints and use injection shift factors ({ISF}s) for network representation. We denote by $\Psi^{n}_{\ell}$ the {ISF} of line $\ell$ with respect to a change in injection at node $n$. 
(See \cite{vanhorn} for the derivation of {ISF}s based on the DC power flow model.) In \eqref{pf}, we represent the real power flow on each line $\ell$ in terms of nodal injections and {ISF}s and constrain it to be within its real power line flow limits. Finally, we impose a nonnegativity constraint on $p_g[k]$, $p^{s}_g[k]$, and $p^{n}_c[k]$ in \eqref{nn1}, \eqref{nn2}, and \eqref{nn3}, respectively. 
\subsection{Evaluation of {FRP} Requirements}\label{3b}
The $\mathsf{SUC}$ problem is an ideal case whose solution yields optimal commitment decisions for {DG}s and optimal decisions on served (i.e., not curtailed) net load levels for each realization. If it were possible to employ this ideal case, $\mathsf{SUC}$ would deliver the minimum expected total cost, under the assumption that it utilizes the true underlying probability distribution of net load. Unfortunately, employing the $\mathsf{SUC}$ solution to clear the market would necessitate radical changes to existing {DAM} practices and so remains unrealistic under current conditions.\par
As such, we treat the optimal $\mathsf{SUC}$ decisions as the gold standard from which we derive the FRP requirements, and we accordingly formulate the $\mathsf{DAMC}$ problem in such a way that brings the optimal $\mathsf{DAMC}$ decisions, as close as possible, toward the optimal $\mathsf{SUC}$ decisions. {W}e set the {FRP} requirements so that the optimal net load levels that are met in the optimal second-stage solution to $\mathsf{SUC}$ can analogously be served by {DG}s that clear the {DAM} under the $\mathsf{DAMC}$ solution. {Note that the optimal $\mathsf{SUC}$ solution for a particular realization could wind up curtailing load, thus only partially meeting a ramping requirement. Such a case would materialize if meeting the ramping requirement in its entirety entailed commitment and generation decisions that markedly drive up the expected total cost. By utilizing the optimal $\mathsf{SUC}$ solution rather than the original net load realizations, we set the FRP requirements so that the ramping requirements are met as long as doing so yields a lower expected total cost.}\par
Specifically, we analyze the change in served net load levels across all consecutive subperiods in all scenarios and determine the maximum and minimum change in net load between each subperiod and its subsequent subperiod over all scenarios. To ensure that the hourly {FRP} requirements for each hour $h \in \mathscr{H}$ can meet the change in the served net load levels within each of the $K$ subperiods, we multiply the maximum and minimum change in net load within each hour by 
{$K$} to set the hourly up-{FRP} and down-{FRP} requirements, respectively. \par
{The FRP requirements are analytically expressed as}
\begin{IEEEeqnarray}{lrCl}
\hspace{-0.5cm}&\underline{\rho}^{\uparrow}(h)  &=&\quad K \cdot\Big( \quad\underset{k \in \mathscr{K}_h}\max\quad\underset{i \in \mathscr{I}}\max\,\sum_{n\in\mathscr{N}} \Big[\big(\xi_{i}^{n}[k+1] \nonumber \\
\hspace{-0.5cm}&&&\hspace{0.05cm}- p^{n}_c[k+1]\Big|_{{y}_{i} = y^{*,\mathsf{SUC}}_{i}}\big) - \big(\xi_{i}^{n}[k] - p^{n}_c[k]\Big|_{{y}_{i} = y^{*,\mathsf{SUC}}_{i}}\big)\Big]\Big),\label{upfrp-1}\\
\hspace{-0.5cm}&\underline{\rho}^{\uparrow}(h)  &=&\quad K \cdot\Big( \,\underset{\begin{subarray}{c} k \in \\  \mathscr{K}_h\setminus\{K\cdot h\}\end{subarray}}\max\;\underset{i \in \mathscr{I}}\max\,\sum_{n\in\mathscr{N}} \Big[\big(\xi_{i}^{n}[k+1] \nonumber \\
\hspace{-0.5cm}&&&\hspace{0.05cm} - p^{n}_c[k+1]\Big|_{{y}_{i} = y^{*,\mathsf{SUC}}_{i}}\big) - \big(\xi_{i}^{n}[k] - p^{n}_c[k]\Big|_{{y}_{i} = y^{*,\mathsf{SUC}}_{i}}\big)\Big]\Big),\label{upfrp-2}\\
\hspace{-0.5cm}&\underline{\rho}^{\downarrow}(h)  &=&-K \cdot\Big( \quad\underset{k \in \mathscr{K}_h}\min\quad\underset{i \in \mathscr{I}}\min\,\sum_{n\in\mathscr{N}} \Big[\big(\xi_{i}^{n}[k+1] \nonumber \\
\hspace{-0.5cm}&&&\hspace{0.05cm}- p^{n}_c[k+1]\Big|_{{y}_{i} = y^{*,\mathsf{SUC}}_{i}}\big) - \big(\xi_{i}^{n}[k] - p^{n}_c[k]\Big|_{{y}_{i} = y^{*,\mathsf{SUC}}_{i}}\big)\Big]\Big),\label{dwfrp-1}\\
\hspace{-0.5cm}&\underline{\rho}^{\downarrow}(h)  &=&-K \cdot\Big( \,\underset{\begin{subarray}{c} k \in \\  \mathscr{K}_h\setminus\{K\cdot h\}\end{subarray}}\min\;\underset{i \in \mathscr{I}}\min\,\sum_{n\in\mathscr{N}} \Big[\big(\xi_{i}^{n}[k+1]  \nonumber\\
\hspace{-0.5cm}&&&\hspace{0.05cm} - p^{n}_c[k+1]\Big|_{{y}_{i} = y^{*,\mathsf{SUC}}_{i}}\big) - \big(\xi_{i}^{n}[k] - p^{n}_c[k]\Big|_{{y}_{i} = y^{*,\mathsf{SUC}}_{i}}\big)\Big]\Big),\label{dwfrp-l}
\end{IEEEeqnarray}
where $x^{*,\mathsf{SUC}}$ denotes the optimal first-stage decisions of the $\mathsf{SUC}$ problem, $y^{*,\mathsf{SUC}}_{i}$ denotes the optimal second-stage decisions of the $\mathsf{SUC}$ problem corresponding to the net load realization $\xi_{i}$, (\ref{upfrp-1}, \ref{dwfrp-1}) hold for all hours $h \in \mathscr{H}\setminus\{24\}$ and (\ref{upfrp-2}, \ref{dwfrp-l}) hold for hour 24.
\subsection{Day-ahead market clearing ($\mathsf{DAMC}$) problem formulation}\label{3c}
We devote this subsection to the description of the $\mathsf{DAMC}$ problem, which is formulated as
\begin{IEEEeqnarray}{l}
\hspace{-0.85cm}\underset{\hspace{0.95cm}\begin{subarray}{l}u_g(h), v_g(h), p_g(h),\\ p_g^{s}(h), p^{n}_c(h), d^{n}(h), \\ r^{\uparrow}_{\wr}(h), r^{\downarrow}_{\wr}(h)\end{subarray}}{\text{minimize}} \sum_{h \in \mathscr{H}}\bigg[ \sum_{g \in \mathscr{G}} \Big[ \alpha^{u}_{g} u_{g}(h)+ \alpha^{v}_{g} v_{g}(h)  + \sum_{s \in \mathscr{S}_g} {\alpha}^{s}_{g}p^{s}_{g}(h)\Big]\nonumber\\
 + \sum_{n \in \mathscr{N}} \alpha_{c}p^{n}_{c}(h)  + \,\alpha_{r} \Big(r^{\uparrow}_{\wr}(h) + r^{\downarrow}_{\wr}(h)\Big)\bigg], \label{damcobj}
\end{IEEEeqnarray}
\begin{IEEEeqnarray}{{-0.01pt}l}
\text{subject to} \nonumber\\
u_{g}(h)-u_{g}(h-1) = v_{g}(h) - w_{g}(h), \label{damcbin}\\
u_{g}(h)-u_{g}^{\circ} = v_{g}(h) - w_{g}(h), \label{damcst}\\
\sum_{h'=h-T^{\uparrow}_{g}+1}^{h} \hspace{-0.5cm}v_{g}(h') \leq u_{g}(h)\hspace{2.7cm} \forall h \in \{T^{\uparrow}_{g},\ldots,24\},\label{damcupt}  \\
\sum_{h'=h-T^{\downarrow}_{g}+1}^{h} \hspace{-0.5cm}w_{g}(h') \leq1- u_{g}(h)\hspace{2.15cm} \forall h \in \{T^{\downarrow}_{g},\ldots,24\}, \label{damcdwt}  \\
\sum_{h'=1}^{\min\{u_g^{\circ}\cdot(T^{\uparrow}_{g} - T^{\uparrow, \circ}_{g}), \,24\}} \hspace{-0.5cm}w_{g}(h') = 0, \label{damciupt}  \\
\sum_{h'=1}^{\min\{(1-u_g^{\circ})\cdot(T^{\downarrow}_{g} - T^{\downarrow, \circ}_{g}), \,24\}} \hspace{-1.0cm}v_{g}(h') = 0, \label{damcidwt}  \\
u_{g}(h), v_{g}(h), w_{g}(h)  \in \{0,1\}  ,\label{damcbinreq}\\
p_{g}(h) \leq (\overline{P}_{g}-\underline{P}_{g})u_{g}(h) ,\label{damcgl} \\
p_{g}(h) \leq p_{g}^{\circ} + {\Delta^{\uparrow}_{g}}u_{g}^{\circ}+(\Delta^{\uparrow,\circ}_{g}-\underline{P}_{g})v_{g}(h),\label{damcrul1} \\
p_{g}(h) \geq p_{g}^{\circ} - {\Delta^{\downarrow}_{g}}u_{g}^{\circ}-(\Delta^{\uparrow,\circ}_{g}-{\Delta^{\uparrow}_{g}}-\underline{P}_{g})w_{g}(h),\label{damcrll1} \\
p_{g}(h) \leq p_{g}(h-1) + {\Delta^{\uparrow}_{g}}u_{g}(h-1)+(\Delta^{\uparrow,\circ}_{g}-\underline{P}_{g})v_{g}(h),\label{damcrul} \\
p_{g}(h) \geq p_{g}(h-1) - {\Delta^{\downarrow}_{g}}u_{g}(h-1)-(\Delta^{\uparrow,\circ}_{g}-{\Delta^{\uparrow}_{g}}-\underline{P}_{g})w_{g}(h),\label{damcrll} \\
p_{g}(h) = \sum_{s \in \mathscr{S}_g} p_{g}^{s}(h) ,\label{damcls} \\
p^{s}_{g}(h) \leq \overline{P}^{s}_{g} - \overline{P}^{\,s-1}_{g}  \hspace{4.25cm} \forall {s} \in \mathscr{S}_g,\label{damclsl} \\
p_g(h) \in \mathbb{R}_{+}, \label{damcnn1}\\
p^{s}_g(h) \in \mathbb{R}_{+}\hspace{5.4cm} \forall {s} \in \mathscr{S}_g, \label{damcnn2}\\
p^{n}(h) = \sum_{g \in \mathscr{G}^{n}} p_{g}(h) + p^{n}_{c}(h) - d^{n}(h) \hspace{2.05cm} \forall n \in \mathscr{N},\label{damcnetpower} \\
\sum_{n \in \mathscr{N}} p^{n}(h)  = 0, \label{damcpb}\\
d^{\,n}(h) = \hat{d}^{\,n}(h)  \quad \xleftrightarrow{\hspace{.5cm}} \quad \lambda^{n}(h) \hspace{2.7cm}  \forall n \in \mathscr{N},  \label{damclmp}\\
\underline{f}_{\ell} \leq \sum_{n \in \mathscr{N}}  \Psi_{\ell}^{n}p^{n}(h) \leq \overline{f}_{\ell} \hspace{3.7cm}  \forall \ell \in \mathscr{L}, \label{damcpf}\\
p^{n}_c(h) \in \mathbb{R}_{+} \hspace{5.5cm} \forall n \in \mathscr{N},\label{damcnn3}\\
-\Delta^{\downarrow}_{g}u_g(h) + (\Delta^{\downarrow}_{g}-\Delta^{\downarrow,\circ}_{g})w_g(h+1) + \underline{P}_g v_{g}[h+1] \leq r_{g}^{\uparrow}(h)\nonumber\\
\hspace{2.85cm}\leq \Delta^{\uparrow}_{g} u_{g}(h+1) +  (\Delta^{\uparrow, \circ}_{g}-\Delta^{\uparrow}_{g})v_g(h+1),\label{upfrplimit}\\
 -\Delta^{\uparrow}_{g}u_g(h+1) + (\Delta^{\uparrow}_{g}-\Delta^{\uparrow,\circ}_{g})v_g(h+1)  \leq r_{g}^{\downarrow}(h)\nonumber \\
 \hspace{1.30cm}\leq  \Delta_{g}^{\downarrow} u_{g}(h) + (\Delta^{\downarrow,\circ}_{g}-\Delta^{\downarrow}_{g})w_g(h+1) - \underline{P}_g v_{g}(h+1) ,\label{dwfrplimit}\\
- \underline{P}_g + \underline{P}_g u_{g}(h+1) \leq r^{\uparrow}_{g}(h) + p_{g}(h) \nonumber\\
 \hspace{2.65cm}\leq  \overline{P}_g - \underline{P}_g u_g(h)+ (\Delta_{g}^{\uparrow,\circ} - \overline{P}_g) v_{g}(h+1), \label{uplimits}\\
- \underline{P}_g + \underline{P}_g u_{g}(h+1) \leq -r^{\downarrow}_{g}(h) + p_{g}(h) \nonumber\\
\hspace{2.65cm}\leq  \overline{P}_g - \underline{P}_g u_g(h)+ (\Delta_{g}^{\uparrow,\circ} - \overline{P}_g) v_{g}(h+1) ,\label{dwlimits}\\
\sum_{g \in \mathscr{G}} r^{\uparrow}_g(h) + r^{\uparrow}_{\wr}(h) \geq  \underline{\rho}^{\uparrow}(h) \quad \xleftrightarrow{\hspace{.5cm}} \quad \varphi^{\uparrow}(h),\label{damcupfrp}\\
\sum_{g \in \mathscr{G}} r^{\downarrow}_g(h) + r^{\downarrow}_{\wr}(h) \geq  \underline{\rho}^{\downarrow}(h) \quad \xleftrightarrow{\hspace{.5cm}} \quad \varphi^{\downarrow}(h),\label{damcdwfrp}\\
u_g(h) \geq u_{g}[h\cdot K]\Big|_{x=x^{*,\mathsf{SUC}}}, \label{sucsols}
\end{IEEEeqnarray}
where \eqref{damcbinreq}--\eqref{damcgl} and \eqref{damcls}--\eqref{sucsols} hold for all hours $h \in \mathscr{H}$, \eqref{damcst} and \eqref{damcrul1}--\eqref{damcrll1} hold for $h \in \{1\}$, and \eqref{damcbin} and \eqref{damcrul}--\eqref{damcrll} hold for $h \in \mathscr{H} \setminus \{1\}$. The objective of the $\mathsf{DAMC}$ problem is to minimize the commitment, startup, and dispatch costs of {DG}s plus the penalty costs due to involuntary load curtailment and {FRP} shortfall. We assume each load has an infinite willingness to pay, and the net demand at each node $\hat{d}^{n},\;\forall n \in \mathscr{N}$, is perfectly inelastic. The variables and constraints introduced in the $\mathsf{DAMC}$  formulation mirror those of $\mathsf{SUC}$ at an hourly granularity. We state in \eqref{damcupt}--\eqref{damcidwt} the minimum uptime and downtime constraints of {DG}s, in \eqref{damcls}-\eqref{damclsl} constraints on the power from each linear segment and the power above minimum, and in \eqref{damcnetpower}-\eqref{damcpf} constraints on net nodal injections and power flows in the network. \par
We bound from above and below the {FRP} provided by each {DG} $g \in \mathscr{G}$ based on its commitment, startup, and shutdown variables in \eqref{upfrplimit} and \eqref{dwfrplimit}, and we jointly constrain the {FRP} provided by each {DG} in each hour $h \in \mathscr{H}$ and its change in power output between each two consecutive hours in \eqref{uplimits}-\eqref{dwlimits}. The equations \eqref{upfrplimit}-\eqref{dwlimits} are based on the three-binary formulation of the {FRP} constraints set forth in \cite{ramp:hobbs}. \par
In our methodology, the up-{FRP} and down-{FRP} requirements are evaluated based on \eqref{upfrp-1}--\eqref{dwfrp-l} as expounded in Section \ref{3b} and are imposed as soft constraints in the $\mathsf{DAMC}$ formulation in \eqref{damcupfrp} and \eqref{damcdwfrp}, respectively. The terms $r^{\uparrow}_{\wr}(h)$ and $r^{\downarrow}_{\wr}(h)$ denote the system-wide up-{FRP} and down-{FRP} shortfall in hour $h \in \mathscr{H}$, respectively.\par
To ensure compatibility with the DAM design principles, we formulate the $\mathsf{DAMC}$ problem at an hourly granularity and employ the hourly net demand $\hat{d}^{\,n}(h)$ in clearing the {DAM}, which is a bid-in quantity that often differs from the intra-hourly net load realizations. In the event that the bid-in net demand is markedly lower than the net load realizations, imposing {FRP} requirements through \eqref{uplimits}--\eqref{dwlimits} alone will not suffice to serve the net load levels in the {RTM} that were met in the optimal $\mathsf{SUC}$ solution. To overcome this shortcoming and drive the optimal {DAM} decisions toward the optimal $\mathsf{SUC}$ solution, we include the constraint \eqref{sucsols} in the $\mathsf{DAMC}$ formulation. By committing the {DG}s in $\mathsf{DAMC}$ that were committed in the optimal $\mathsf{SUC}$ solution, \eqref{uplimits}--\eqref{dwlimits} and \eqref{sucsols} jointly ensure that the {DG}s receiving {DAM} awards can serve the net load levels that were met in the optimal $\mathsf{SUC}$ solution for each considered net load realization $\xi^{n}_{i}[k],\,i\in\mathscr{I}$. As such, albeit the use of the bid-in net demand in the second pass, our methodology can bring forth {DAM} awards that help bring operating costs closer to the minimum expected total cost of the $\mathsf{SUC}$ solution.\par
In line with the current {ISO}/{RTO} market models, the objective function \eqref{damcobj} does not include a term for the dispatch costs that may arise due to the deployment of procured {FRP}s. As such, if we were to enforce {FRP} requirements through \eqref{uplimits}--\eqref{dwlimits} without stipulating \eqref{sucsols}, we could potentially observe undesirable {DAM} schedules. Specifically, in bringing online additional {DG}s to meet the {FRP} requirements, the optimal {DAM} solution could prioritize {DG}s with lower startup and commitment costs without taking into account their dispatch costs. Such a {DAM} solution would clearly run counter to our objective of driving the {DAM} decisions toward the optimal $\mathsf{SUC}$ solution; furthermore, it could bring substantial dispatch costs that would have to be borne if the acquired {FRP}s were called upon. The constraint \eqref{sucsols} aids our methodology in surmounting these challenges, as it enforces that the $\mathsf{DAMC}$ solution prioritize for {DAM} awards those {DG}s that were committed under the optimal solution to the $\mathsf{SUC}$ problem, which intrinsically assesses the dispatch costs of {DG}s in conjunction with the commitment and startup costs. Owing to \eqref{sucsols}, $\mathsf{DAMC}$ can capitalize on this cost assessment in obviating unexpectedly high dispatch costs, all the while tapering off the total operating cost of $\mathsf{DAMC}$ closer to that of $\mathsf{SUC}$. \par
Our methodology further takes into account the nodal deliverability of the procured {FRP}s.
{Since the $\mathsf{SUC}$ problem explicitly includes power flow constraints for all scenarios and intra-hourly subperiods}, our methodology ensures that the {DG}s that were committed in the $\mathsf{SUC}$ solution and receive {FRP} awards under the $\mathsf{DAMC}$ solution can deliver their {FRP}s to serve the net load levels that were met in the $\mathsf{SUC}$ solution.\par 
We model the nodal net demand ${d}^{\,n}(h)$ to be a continuous variable and, in \eqref{damclmp}, set it equal to the parameter $\hat{d}^{\,n}(h)$, which denotes the {bid-in} net demand at node $n$ in hour $h$. The shadow price\footnote{The $\mathsf{DAMC}$ problem is formulated as a mixed-integer linear problem, for which dual variables are not well defined. To evaluate the shadow prices, we follow a common practice \cite{oneill} also used by most North American {ISO}s/RTOs, whereby we first solve the $\mathsf{DAMC}$ problem, fix the binary variables to their optimal values, relax the binary variables to being continuous, solve the resulting continuous convex problem, and finally evaluate the shadow prices as a by-product of the continuous problem solution.} for the constraint \eqref{damclmp} for node $n$ and hour $h$, expressed by $\lambda^{n}(h)$, denotes the sensitivity of the cost function with respect to the value of net demand at node $n$ and hour $h$ and as such directly corresponds to the {LMP} at node $n$ for hour $h$ \cite{m39rep}. We further utilize the shadow price for the constraints \eqref{damcupfrp} and \eqref{damcdwfrp} to set the price for up-{FRP} and down-{FRP} in hour $h$, respectively. The up-{FRP} (resp. down-{FRP}) price represents the opportunity cost of dispatching downward a cheap {DG} and in turn dispatching upward an expensive {DG} so as to meet up-{FRP} (resp. down-{FRP}) requirements, and is equal to the difference between the dispatch cost of the cheap {DG} and 
{energy price}.\footnote{See \cite[Sec. 6]{ramp:add} and Footnote 6 in \cite{ramp:hobbs} for a detailed discussion on the implications of {FRP} prices.} We do not include spinning reserves, virtual bidding, or contingencies in our formulation, as their inclusion would complicate the $\mathsf{DAMC}$ formulation without shedding light on the fundamental issues concerning {FRP}s.
\subsection{Computational Issues}
{The capability of the proposed methodology to assist SOs in real-life applications is contingent upon its computation time. While the development of specific methods that target the efficient application of our methodology for real-life systems is beyond the scope of this article, we point out studies from the literature that could prove useful to this end.}\par
{Utilizing a greater number of scenarios in our methodology could help more effectively capture the uncertainty of net load and thus lead to more judicious FRP requirements. Nevertheless, under increasing number of scenarios, solving the $\mathsf{SUC}$ problem could be computationally prohibitive and thus hamper the feasibility of our methodology. To address this problem, }\cite{romisch14} {proposes scenario reduction techniques to represent uncertainty with fewer scenarios. Since the performance of a stochastic programming model greatly depends on how well the utilized scenarios represent the underlying uncertainty,} \cite{ershun} {puts forward a scenario mapping technology that aims to reduce model complexity---all the while preserving the uncertainty of wind generation to the greatest extent. Studies in} \cite{knn, powertech} {utilize the k-means clustering algorithm to partition similar scenarios into clusters, among which} \cite{ powertech} {employs distance measures tailored to quantify the similarity between time series. Such techniques could be leveraged to determine those few net load scenarios whose utilization in solving the $\mathsf{SUC}$ problem ultimately yields FRP procurement levels that can meet both the variability and uncertainty in net load. In addition, a great deal of studies in the literature works out decomposition techniques for efficiently applying stochastic programming models to take UC decisions in real-life systems, such as decomposing the problem on the basis of units} \cite{unit_dec}{, scenarios} \cite{scen_dec}{, and time periods }\cite{temp_dec}{. Drawing upon these techniques in applying our methodology will allow for a faster solution of $\mathsf{SUC}$, which will greatly aid in the application of our approach for real-life systems.} 
\section{Numerical Experiments}\label{4}
In this section, we carry out numerical studies using various benchmarks and evaluation metrics to demonstrate the application and effectiveness of the proposed methodology.
\subsection{Benchmark Methods}\label{4a}
To illustrate the relative merits of the proposed methodology, we include as benchmark methods three different approaches to {FRP} procurement. The benchmark method termed $\mathsf{DAMC-nf}$ ($\mathsf{nf}$ as in not fixed) procures {FRP}s in the {DAM} based on the levels evaluated by \eqref{upfrp-1}--\eqref{dwfrp-l}, yet it does not fix the commitment status of the {DG}s based on the optimal $\mathsf{SUC}$ decisions and as such does not include \eqref{sucsols} in its formulation. The formulation for the $\mathsf{DAMC-nf}$ problem is expressed as
\begin{IEEEeqnarray}{ll}
\hspace{-5.4cm}\text{minimize} &\hspace{-3.5cm} \eqref{damcobj}\nonumber\\
\hspace{-5.4cm}\text{subject to} & \hspace{-3.5cm} \eqref{damcbin}\text{--}\eqref{damcdwfrp}\nonumber\\
& \hspace{-3.5cm} \eqref{upfrp-1}\text{--}\eqref{dwfrp-l}.\nonumber
\end{IEEEeqnarray}
Another benchmark method we utilize is the $\mathsf{DAMC-95}$ method, which sets the hourly up-{FRP} and down-{FRP} requirement based on the upper and lower bound of the 95\% confidence interval of net load in each hour and procures {FRP}s without considering the commitment decisions taken under the $\mathsf{SUC}$ problem. The mathematical formulation of the $\mathsf{DAMC-95}$ method is 
\begin{IEEEeqnarray}{lll}
\hspace{-5.3cm}\text{minimize} &\hspace{0.4cm}& \eqref{damcobj}\nonumber\\
\hspace{-5.3cm}\text{subject to}&& \eqref{damcbin}\text{--}\eqref{damcdwfrp}.\nonumber
\end{IEEEeqnarray}
To insulate the effects of acquired {FRP} capacity on the results, we further utilize a benchmark method termed $\mathsf{DAMC-w/o}$ that does not procure {FRP}s, expressed as
\begin{IEEEeqnarray}{lll}
\underset{\begin{subarray}{c}u_g(h), v_g(h),\\ p_g(h), p_g^{s}(h),\\ p^{n}_c(h)\end{subarray}}{\text{minimize}} &\hspace{0.4cm}& \sum_{h \in \mathscr{H}}\bigg[ \sum_{g \in \mathscr{G}} \Big[ \alpha^{u}_{g} u_{g}(h) + \alpha^{v}_{g} v_{g}(h) + \sum_{s \in \mathscr{S}_g} {\alpha}^{s}_{g}p^{s}_{g}(h)\Big]   \nonumber\\[5pt]
&\hspace{0.4cm}&+ \sum_{n \in \mathscr{N}} \alpha_{c}p^{n}_{c}(h)\bigg] \label{damcwoobj}\\[5pt]
\text{subject to} &\hspace{0.4cm}& \eqref{damcbin}\text{--}\eqref{damcnn3}.\nonumber
\end{IEEEeqnarray}
\subsection{Evaluation Metrics}\label{4b}
In this subsection, we introduce evaluation metrics that we draw upon to compare the performance of the proposed $\mathsf{DAMC}$ methodology with that of the benchmark methods set forth in Section \ref{4a}. In our numerical studies, we employ a two-settlement system that is standard in U.S. organized markets. The solution of the $\mathsf{DAMC}$ problem as well as that of each benchmark method constitutes a forward {DAM} with financially binding awards, and the deviations from the day-ahead schedules of each method are exposed to the prices of its associated {RTM}. While {RTM}s in different jurisdictions may involve various processes, we solely consider the fifteen-minute market ({FMM}). We assume that the commitment decisions of the {DAM} cannot be altered in the {FMM}. \par
For simplicity, we model the {FMM} with a two-hour scheduling horizon and eight subperiods, and we run the {FMM} at the top of each hour, thus running a sequence of 24 {FMM}s for a given day. We further suppose that, in each {FMM} run, the decisions for only the first four subperiods are binding. We do not procure {FRP}s in the {FMM} so as to isolate the impact of {FRP} procurement in the day-ahead time frame. Each simulation includes the same set of {DG}s for the {DAM} and its {FMM} in order to highlight the performance of the methods themselves and preclude the results from being affected by the characteristics of fast-ramping resources that turn on in the {FMM} after the {DAM} results are published. Commensurate with the choice of {FMM}, we consider four intra-hourly subperiods in the $\mathsf{SUC}$ formulation, that is, $K=4$.  \par
In our experiments, we model the net load at each node for each intra-hourly subperiod to be normally distributed with mean $\overline{\xi}_{n}[k]$ and standard deviation equal to one percent of its mean, i.e.,  ${\xi}_{n}[k] \sim \mathcal{N}\big(\overline{\xi}_{n}[k],\,\overline{\xi}_{n}[k]\cdot1\%\big)$. We further model the hourly net demand at each node for the $\mathsf{DAMC}$ method, as well as the benchmark methods, to be equal to the statistical average of the mean value of net load in its $K$ constituent intra-hourly subperiods, i.e, $\hat{d}^{n}(h) = \frac{1}{K} \sum_{k \in \mathscr{K}_h} \overline{\xi}_{n}[k],\;\forall h \in \mathscr{H}$. It is worth noting that the methodology presented in Section \ref{3} is generic, in that it can be tailored to any choice of representation of the uncertain intra-hourly net load ${\xi}[k]$ as well as the hourly net demand $\hat{d}^{n}(h)$. To evaluate the out-of-sample performances, we draw a new \textit{i.i.d.} random sample of $R'$ realizations of the random net load $\xi$, generated independently of the sample used to solve the $\mathsf{SUC}$ problem. For each method and each of the $R'$ realizations, we clear the 24 {FMM}s and record the results for the binding intervals.\par
A metric that is commonly used to assess the performance of market-based solutions is social welfare, i.e., producer surplus plus consumer surplus. Nevertheless, social welfare is not well defined in our case as the demand is modeled to be inelastic and so the consumer surplus of all methods goes to infinity. As such, in line with previous work \cite{transdereg} that modeled the demand to have zero price elasticity, we take the social welfare as the negative of the total payments, and compute the total payments in all binding subperiods of the 24 {FMM} runs over all $R'$ realizations. Note that the total payments include uplifts in the form of make-whole payments, whereby each {DG} that incurs costs higher than the payments it receives for its energy and {FRP} awards receives an uplift payment so that it breaks even. In addition, we employ as an evaluation metric the total curtailed load in all binding subperiods of the 24 {FMM} runs over all $R'$ realizations. 
\subsection{Results}\label{4c}
We tested our methodology on several systems that varied in size from 1 to 118 buses. To demonstrate the effectiveness of our methodology in a clear manner, we present in this article the results on a modified IEEE 14-bus test system \cite{matpower} with 5 {DG}s, 20 transmission lines, and 11 loads\footnote{The simulation data is provided in \url{https://github.com/oyurdakul/damc}.} which has sufficient scope to illustrate the capabilities of our methodology. The results are representative and are not limited to systems of this size. We solve the $\mathsf{SUC}$ problem for $R=100$ and conduct the RTM out-of-sample assessments for $R’=50$. We build all models in Julia and modify the \verb|UC.jl| package \cite{ucjulia} for our tests by extending it to the two-stage stochastic setting and incorporating {FRP}s to its base formulation. We use Gurobi 9.5 for solving all optimization problems and conduct the simulations on a 2.6 GHz Intel Core i7 CPU with 16 GB of RAM.\par
\begin{table}[h]
\centering

\setlength{\tabcolsep}{7pt} 
\renewcommand{\arraystretch}{1.5}
\caption{Case study results}
\label{results}
\centering
\begin{tabular}{ c | c | c }
\hline \hline 
& \multirow{2}{*}{total payment} &  total load\\[-5pt]
&& curtailment ({MW})\\\hline\hline
 $\mathsf{DAMC}$ &  \$19,942,069.96 & 2.21 \\\hline
 $\mathsf{DAMC-nf}$ & \$20,383,836.90 & 9.17 \\\hline
   $\mathsf{DAMC-95}$ & \$24,483,164.43 & 376.58\\\hline
 $\mathsf{DAMC-w/o}$ & \$27,962,586.92 & 766.71\\\hline
\hline \hline
\end{tabular}
\end{table}
We provide in Table \ref{results} the results under the proposed $\mathsf{DAMC}$ formulation as well as the benchmark methods for different evaluation metrics. Table \ref{results} makes clear the monetary benefits that can be reaped by applying our methodology. The proposed methodology delivers the lowest total payment and notches a 2.17\% reduction in total payment compared to the runner-up $\mathsf{DAMC-nf}$ method. These results bear out the central role of taking into account the dispatch costs that would be incurred if the procured {FRP}s were to be deployed. The lower total payment achieved by the $\mathsf{DAMC-nf}$ method vis-à-vis the other benchmark methods drives home the merits of evaluating the {FRP} requirements as per the results of the $\mathsf{SUC}$ problem in lieu of pre-specified limits, such as the 95\% confidence interval. Among all methods, the $\mathsf{DAMC-w/o}$ method yields the highest total payment, which brings out the importance of {FRP} procurement in the {DAM}. The higher performance of our methodology based on this monetary metric is echoed in fending off involuntary load curtailment, as the total curtailed load attains its lowest value under the proposed methodology, followed by the $\mathsf{DAMC-nf}$ method, yet it climbs sharply under the $\mathsf{DAMC-w/o}$ method. \par
\begin{figure}[h]
\includegraphics[width=.95\linewidth]{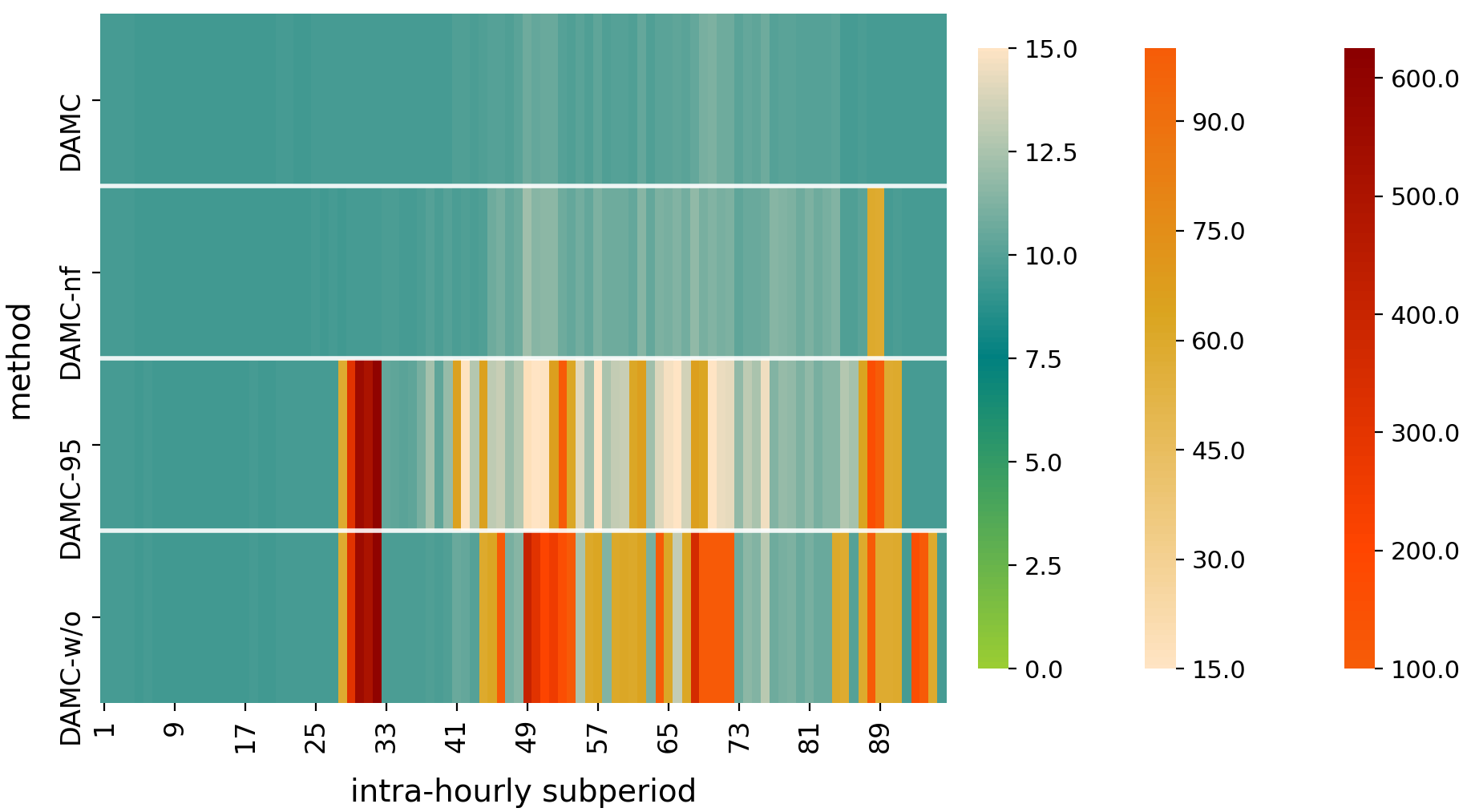}
\caption{The average {LMP} of all system buses in the {FMM}}
\label{lmps}
\end{figure}
We next investigate the {LMP}s over the simulation horizon for the {FMM} of each method. In Fig. \ref{lmps}, we illustrate the average {LMP} of all system buses for the {FMM} of each method in each 15-minute subperiod. We observe in Fig. \ref{lmps} that the {LMP}s under the proposed methodology are, by and large, lower than those under the benchmark methods. Most notably, the {LMP}s under the benchmark methods surge in the subperiods of hours 8, 13, and 23, whereas those under the proposed methodology remain at a relatively low ebb. \par
\begin{figure}[h]
\hspace{-0.1cm}
\includegraphics[width=.95\linewidth]{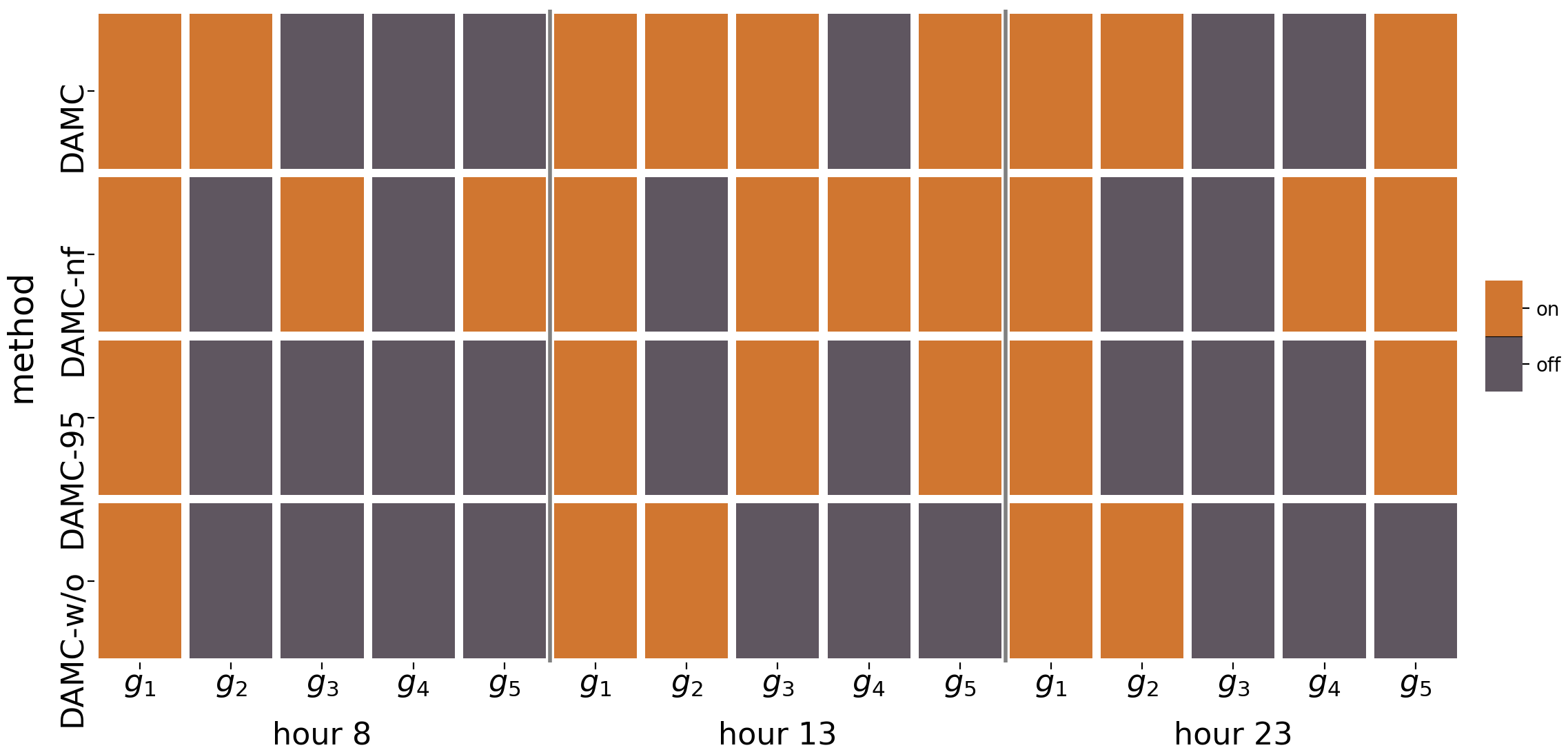}
\caption{The commitment decisions in hours 8, 13, and 23}
\label{coms}
\end{figure}
To gain insights into the factors that inflate the {LMP}s under these methods, we depict in Fig. \ref{coms} the commitment status of the five {DG}s under all methods in hours 8, 13, and 23. Fig. \ref{coms} shows that, in hour 8, $g_1$ is committed under all methods. Additionally, $\mathsf{DAMC}$ commits $g_2$, whereas $\mathsf{DAMC-nf}$ commits $g_3$ and $g_5$, and $\mathsf{DAMC-95}$ and $\mathsf{DAMC-w/o}$ commit no additional {DG}s. In a similar vein, in hour 13, $g_1$, $g_3$, and $g_5$ are committed by all methods except $\mathsf{DAMC-w/o}$. Additionally, $\mathsf{DAMC}$ commits $g_2$, while $\mathsf{DAMC-nf}$ commits $g_4$, and $\mathsf{DAMC-95}$ does not commit further DGs. In all these three hours, $\mathsf{DAMC-w/o}$ commits the fewest DGs. \par
We ascribe the prioritized commitment of $g_2$ in lieu of $g_3$, $g_4$, and $g_5$ under $\mathsf{DAMC}$ to $\mathsf{DAMC}$'s capability in effectively exploiting the optimal $\mathsf{SUC}$ decisions. Compared with the other generators, $g_2$ has higher startup and commitment costs yet a lower dispatch cost. On account of leveraging the optimal solution to the $\mathsf{SUC}$ problem, which comprehensively evaluates the commitment, startup, and dispatch costs, $\mathsf{DAMC}$ prudently elects to meet {FRP} requirements through $g_2$. The $\mathsf{DAMC-nf}$ and $\mathsf{DAMC-95}$ methods, however, solely consider the startup and commitment costs in committing {DG}s to meet {FRP} requirements. As such, in hours 8, 13, and 23, $\mathsf{DAMC-nf}$ commits $g_5$, which has a relatively low startup and commitment cost but a high dispatch cost. Ultimately, the preferential commitment of $g_2$ under $\mathsf{DAMC}$ brings forth a lower marginal cost of energy and in turn lower {LMP}s. In light of both the high ramping capability and low dispatch cost of $g_2$ compared with the other DGs, these results attest to the capacity of $\mathsf{DAMC}$ to spur investments in {DG}s that can meet steep ramping requirements at a low cost of dispatch. Furthermore, we point out that the commitment decisions under $\mathsf{DAMC-95}$ and $\mathsf{DAMC-w/o}$ lead to numerous subperiods of involuntary load curtailment and hence to markedly high {LMP}s in the {FMM}. \par
\begin{figure}[h]
\centering
\includegraphics[width=.7\linewidth]{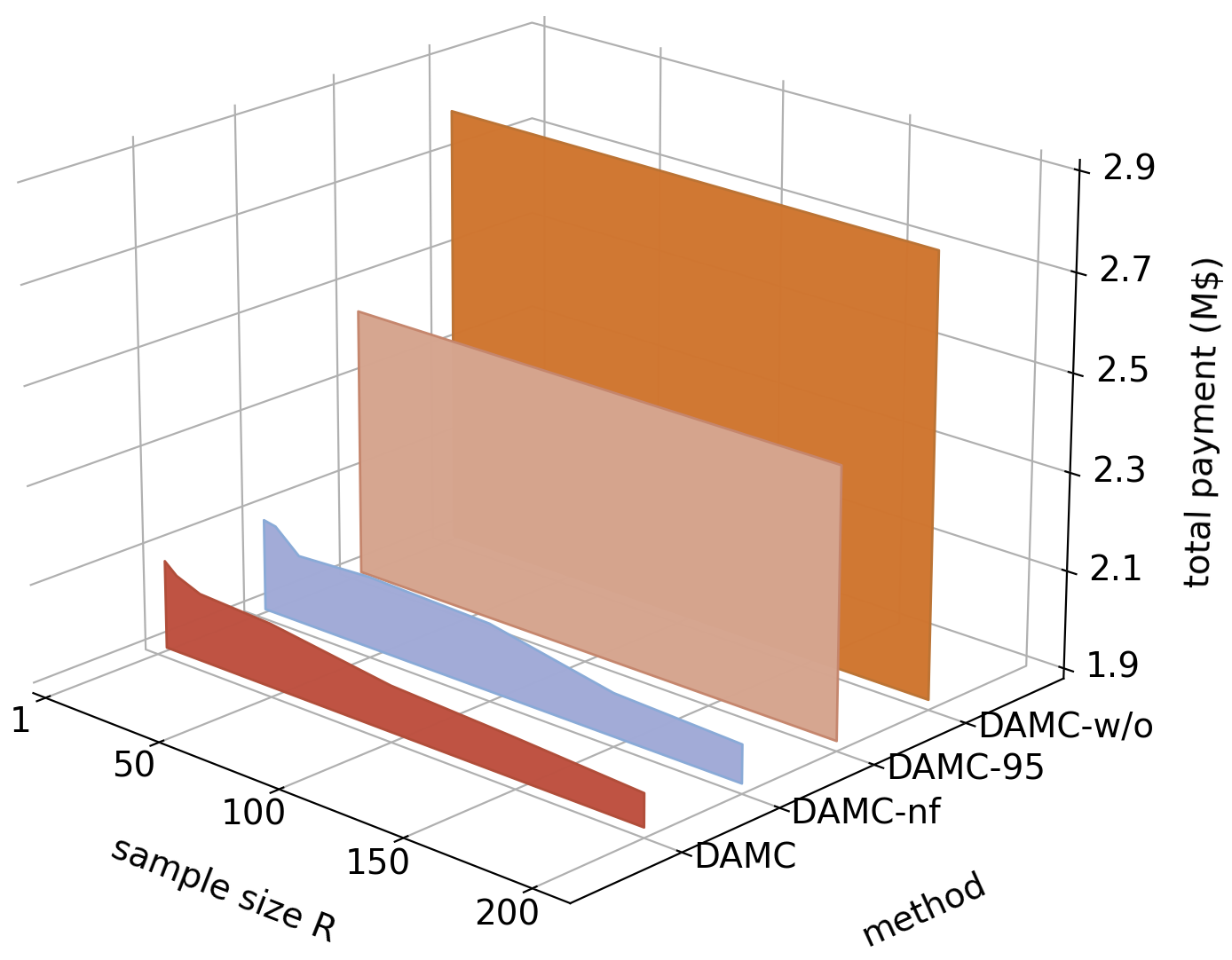}
\caption{Total payment under different methods as a function of sample size}
\label{sens}
\end{figure}
Another assessment of particular interest is whether increasing the size of the \textit{i.i.d.} random sample used in solving the $\mathsf{SUC}$ problem can help the $\mathsf{SUC}$ problem better grasp the underlying uncertainty and variability of net load, thereby further reducing the total payment under the proposed methodology. To this end, we repeat our experiments for the following values of $R=5,10, 20, 50, 100, 150, \text{and }200$. We observe in Fig. \ref{sens} that the total payment under the proposed methodology generally decreases as $R$ increases. This shows that enlarging the sample size augments the ability of $\mathsf{SUC}$ to evaluate the {FRP} requirements and commitment decisions, which is in turn leveraged by the $\mathsf{DAMC}$ methodology. \par
Similar to $\mathsf{DAMC}$, the total payment under $\mathsf{DAMC-nf}$ generally decreases with increasing sample size, yet it exhibits larger fluctuations compared with $\mathsf{DAMC}$. We attribute these larger fluctuations to the fact that the deployment of {FRP}s procured under $\mathsf{DAMC-nf}$ may prompt hefty dispatch costs. By contrast, the total payments under $\mathsf{DAMC-95}$ and $\mathsf{DAMC-w/o}$ remain constant as $R$ increases, since these methods do not take into account the optimal $\mathsf{SUC}$ decisions.  
\section{Conclusion}\label{5}
This paper maps out a methodology for the acquisition of flexible ramping product ({FRP}) in the day-ahead market ({DAM}). The proposed methodology comprises two market passes that are consecutively executed, which allows for not only a comprehensive evaluation of the expected total operating cost, but also a close adherence to today's {DAM} designs. \par
A principal strength of our methodology is the derivation of the {FRP} requirements based on the optimal solution to a stochastic unit commitment ($\mathsf{SUC}$) problem that considers the uncertainty and the inter-hourly and intra-hourly variability of net load. As such, our methodology warrants an economic basis for the stipulated {FRP} levels, which indeed are imposed in a way that pushes the {DAM} awards toward the optimal decision of the $\mathsf{SUC}$ problem, thereby driving the costs closer to the minimum expected total cost obtained under $\mathsf{SUC}$. By assessing the expected dispatch costs in determining {FRP} awards, our approach can forestall extortionate costs that could result from the deployment of procured {FRP}s. We demonstrate the application and effectiveness of the proposed methodology using representative studies. The results lay out the quantifiable improvement that the proposed methodology brings over selected benchmark methods under various evaluation metrics.  \par
A natural extension of this work is to set the requirements of and price {FRP}s on a nodal basis. Another avenue for future research is the utilization of a distributionally robust optimization model in the first market pass to help avoid a poor out-of-sample performance that a stochastic optimization model may exhibit, if the assumed probability distribution for net load is markedly different than its true distribution. {We further set out to develop scenario reduction and decomposition techniques that ensure solving $\mathsf{SUC}$ in an efficient manner, thus providing for the application of our methodology in real-life systems.} 
\bibliographystyle{IEEEtran}
\bibliography{IEEEabrv,rampbib}
%







\end{document}